\documentclass{article}
\usepackage{euscript,amsfonts,amssymb,amsmath}
\usepackage{graphicx,epsfig,amscd}
\usepackage[all]{xy}
\author{A. G. Gorinov}
\title{On the cobordisms of M\"obius circles}
\date{}

\newcommand{\N}{{\mathbb N}}
\newcommand{\Z}{{\mathbb Z}}

\newcommand{\C}{{\mathbb C}}

\newcommand{\R}{{\mathbb R}}

\newcommand{\eu}{\EuScript}

\newtheorem{theorem}{Theorem}
\newtheorem{lemma}{Lemma}
\begin{document}
\maketitle
\begin{abstract}
The boundary of a M\"obius manifold carries a canonical M\"obius structure. This enables one to define the cobordism group of $n$-dimensional (closed)
M\"obius manifolds. The purpose of this note is to show that the cobordism group of M\"obius circles is zero, i.e., every M\"obius circle bounds a M\"obius surface. We also complete N. Kuiper's classification of projective
structures on $S^1$ (we show that there are in fact two series of projective circles with parabolic holonomy, and not one).
\end{abstract}
\section{Introduction}
Recall that a {\it M\"obius} manifold is an oriented $(G,X)$-manifold (see e.g., \cite[3.1.2]{thurston}) of dimension $n$ with $X=S^n$ and $G$ the group $\mbox{M{\"o}b}_n$ of orientation-preserving M\"obius transformations (see e.g., \cite{hert} for details). We suppose that the charts participating in the definition of the M\"obius structure are orientation preserving. Two M\"obius structures on the same manifold are said to be {\it equivalent}
if one of them is induced from the other under some diffeomorphism homotopic to the identity.

For a M\"obius manifold $M$, we set $-M$ to be the M\"obius manifold ``inverse to $M$'', i.e., as a smooth oriented manifold, $-M$ is equal to $M$ with the orientation reversed, and the charts that define the M\"obius structure on $-M$ are the compositions
of the charts for $M$ and an orientation-reversing M\"obius transformation. In the sequel, all M\"obius diffeomorphisms are assumed to be orientation preserving. A diffeomorphism $M_1\to M_2$ between two M\"obius manifolds is said to be {\it anti-M\"obius}, if the composition
$-M_1\stackrel{\mathrm{Id}}\longrightarrow M_1\to M_2$ is M\"obius.

The definition of M\"obius structures can be naturally extended to manifolds with boundary.
Let $M$ is a M{\"o}bius manifold, and let $p:N_1\to N_2$ is an anti-M{\"o}bius diffeomorphism between two disjoint components of $\partial M$. Then the 
manifold $M'$ obtained from $M$ by identifying $x$ and $p(x)$ for all 
$x\in N_1$ can be equipped with a canonical M{\"o}bius structure.

Hence, one can define the {\it cobordism group $\Omega_n^{\mbox{\scriptsize \rm M{\"o}b}}$ of $n$-dimensional M\"obius manifolds} as follows.
We set the elements of $\Omega_n^{\mbox{\scriptsize \rm M{\"o}b}}$ to be $n$-dimensional closed M\"obius manifolds considered modulo
the equivalence relation $\sim$ with $M_1\sim M_2$, iff there exists an $n+1$-dimensional compact M\"obius manifold $N$ such that
$\partial N$ is M\"obius diffeomorphic to $M_1\sqcup -M_2$. The group operation
is induced by taking the disjoint union. The main result of this note is the following theorem.

\begin{theorem}
\label{nqs}$\Omega_1^{\mbox{\scriptsize \rm M{\"o}b}}=0$.
\end{theorem}

The proof is straightforward: we describe in section \ref{sec2} the M\"obius structures on the circle, and for each of those we construct in section \ref{sec3} a bounding surface. In section \ref{sec2}, we also show (lemma \ref{projective} and remark 2) that some cases are missing from N. Kuiper's classification \cite{kuiper} of projective structures
on the circle, and we complete that classification.

{\bf Remark 1.} One can define for any $n\geq 1$ the cobordism category $\eu{MC}_n$ of $n$-dimensional M\"obius manifolds as follows. We set the objects of $\eu{MC}_n$ to be compact M\"obius manifolds of dimension $n$,
and for any two objects $N_1$ and $N_2$ of $\eu{MC}_n$, morphisms from $N_1$ to $N_2$
are defined to be $n+1$-dimensional compact M\"obius manifolds $N$ such that $\partial N=N_1\sqcup -N_2$ (modulo M\"obius diffeomorphisms that are identical on $\partial N$). The categories $\eu{MC}_n$ are strict symmetric monoidal, and
the above theorem \ref{nqs} can be reformulated as $K_0(\eu{MC}_1)=0$. It would be interesting to learn more about $\eu{MC}_1$, e.g., a presentation by generators and relations etc.

The author is grateful to P. Vogel for proposing the problem and useful discussions.
\section{M{\"o}bius structures on the oriented circle}
\label{sec2}
In this section we recall the classification of M\"obius structures on the circle $S^1$.
A M{\"o}bius structure on a 1-dimensional manifold can
be viewed as a $(G,X)$-structure, where 
$X$ is any circle $C$ in $\C P^1$,
and $G$ is the group of automorphisms of $\C P^1$ that
preserve both disks bounded by $C$. We shall often use the projective model, i.e., we shall take the circle $\R P^1\subset\C P^1$, so that $G$ becomes the group of orientation-preserving projective automorphisms of $\R P^1$.
Given a M\"obius structure $S$ on $S^1$, we denote by $[S]$ the equivalence class of $S$, and we define $\eu{M}(S^1)$ to be the set of all equivalence classes of M\"obius structures on $S^1$.

Let us identify the universal cover of $S^1$ with the (canonically oriented) $\R$.
Any M\"obius structure $S$ on $S^1$ constructed below will be defined by
an orientation-preserving developing map
$ F_S$ from $\R$ to a model circle $C\subset\C P^1$ and a holonomy homomorphism $f_S$ from $\Z$ to the subgroup of $\mathrm{PSL}_2(\C)$ corresponding to $C$;
the maps $F_S$ and $f_S$ should satisfy the condition
\begin{equation}
\label{eqqq}
 F_S(x+1)= f_S(1)\cdot F_S(x).
\end{equation}

Let $A_S$ be a matrix
that represents $ f_S(1)$ (this matrix is defined up to a
sign). We 
shall distinguish the following cases:

\begin{enumerate}
\item $|tr(A_S)|\leq 2$, and $A_S$ is diagonalisable.
\item $|tr(A_S)|=2$, $A_S$ is nondiagonalisable.
\item $|tr(A_S)|>2$.
\end{enumerate}

The structures that verify (1), (2), (3), will be called respectively 
{\it elliptic, parabolic} and {\it hyperbolic}.

For any $\alpha>0$ and $n>0$ let $\psi_{n,\alpha}$ be a continuous map $[0,1]\to\R P^1$ such that 1. $\psi_{n,\alpha}=\mathrm{e}^{\alpha x}-1,$ when $x$ is close to 0,
$\psi_{n,\alpha}=\mathrm{e}^{\alpha}(\mathrm{e}^{\alpha x}-1),$ when $x$ is close to 1, 2. $\psi^{-1}_{n,\alpha}(0)$ consists of $n+1$ elements, and 3.
$\psi_{n,\alpha}|(0,1)$ is an orientation-preserving local diffeomorphism.
Analogously, for any couple $(n,\varepsilon)$ (where $n>0$ is an integer, and $\varepsilon=\pm 1$)
let $\xi_{n,\varepsilon}$ be a continuous map $[0,1]\to\R P^1$ such that 1. $\xi_{n,\varepsilon}=-\frac1x,$ when $x$ is close to 0,
$\xi_{n,\varepsilon}=\varepsilon+\frac{1}{1-x},$ when $x$ is close to 1, 2. $\xi^{-1}_{n,\varepsilon}(\infty)$ consists of $n+1$ elements, and 3.
$\xi_{n,\varepsilon}|(0,1)$ is an orientation-preserving local diffeomorphism.
The following theorem gives an explicit representative for each element of $\eu M (S^1)$ (cf. \cite{kuiper}).

\begin{theorem}
\label{circle}
\begin{enumerate}
\item
If two M{\"o}bius structures on $S^1$ belong to different types, they are
not equivalent.
\item Equivalence classes of hyperbolic M{\"o}bius structures are
parametrised by couples $(n,\alpha)$, where $\alpha$ is a positive real number, and $n$ is a nonnegative integer.

The corresponding maps $ F_S:\R\to \R P^1,  f_S:\Z\to \mathrm{PSL}_2(\R)$ can be chosen as follows. Let $f_S(1)$ be the map $z\mapsto \mathrm{e}^{\alpha}z$.
If $n=0$, set $ F_S(x)=\mathrm{e}^{\alpha x}$, otherwise set $ F_S(x)=\psi_{n,\alpha} (x)$ for $x\in [0,1]$ and extend $ F_S$ to $\R$
using (\ref{eqqq}). Denote the equivalence class of hyperbolic M\"obius structures
corresponding to the couple $(n,\alpha)$ by
$H_{n,\alpha}$.
\item
Equivalence classes of parabolic M{\"o}bius
structures are parametrised by elements of the set
$(\{\mbox{nonnegative integers}\}\times\{\pm 1\})\setminus\{(0,-1)\}$.

The equivalence class that corresponds to the couple $(0,1)$ is represented by the M\"obius structure defined by the following maps
$ F_S:\R\to \R P^1,  f_S:\Z\to \mathrm{PSL}_2(\R)$: set $ F_S(x)=x,  f_S(1)(z)=z+1$. 

The equivalence class corresponding to the couple $(n,\varepsilon)$ (where $n$ is a positive integer, and $\varepsilon$ is 1 or $-1$)
is represented by the M\"obius structure defined by the maps $ F_S:\R\to \R P^1,  F_S:\Z\to \mathrm{PSL}_2(\R)$ that
can be constructed as follows: set $ f_S(1)(z)=z+\varepsilon$, set $ F_S(x)=\xi_{n,\varepsilon} (x)$ for $x\in [0,1]$ and extend $ F_S$ to $\R$
using (\ref{eqqq}). Denote the equivalence class of parabolic M\"obius structures
corresponding to the couple $(n,\varepsilon)$ by
$P_{n,\varepsilon}$.
\item Equivalence classes of elliptic M{\"o}bius structures are 
parametrised by positive real numbers. The equivalence class corresponding to $\alpha>0$ can be represented by the M\"obius structure defined by the maps $F_S:\R\to\{\mbox{complex numbers of modulus 1}\}$ and $f_S(1):\Z\to\mathrm{PU}_{1,1}$ given by $F_S(x)=\mathrm{e}^{i\alpha x},f_S(1)(z)=\mathrm{e}^{i\alpha}z$. 
Denote this equivalence class by $E_\alpha$.
\end{enumerate}
\end{theorem}

{\bf Proof of theorem \ref{circle}.}
The argument is analogous to the one used in \cite{kuiper} to classify projective circles modulo
projective diffeomorphisms (but see remark 2 below).

First, let us show that any M\"obius structure $S$ on $S^1$ is equivalent to some structure described in theorem \ref{circle}.
Suppose, e.g., that $S$ is hyperbolic. Choose a base point in $S^1$, and denote by
$ F_S$ the corresponding developing map $\R\to\R P^1$; let $ f_S:\Z\to \mathrm{PSL}_2(\R)$ be the homomorphism
such that $ F_S$ and $ f_S$ satisfy (\ref{eqqq}). Then, replacing $S$ by an equivalent structure, we can assume that
$ f_S$ is given by the formula
$ f_S(n)(z)=\mathrm{e}^{\alpha n}z$ for some $\alpha>0$. The image of $ F_S$ is either the whole $\R P^1$ or not; we can assume (changing the base
point, if necessary) that $ F_S(0)$ is 1 in the first case and 0 in the second.
Now it is easy to see that $S$ is equivalent to some structure from the second
assertion of theorem \ref{circle}. The parabolic and elliptic cases are considered in an analogous way.

Now, let us prove that the structures introduced in theorem \ref{circle} are pairwise nonequivalent. Let $S$ be a M\"obius structure on $S^1$.
Suppose that $S$ is
defined by an orientation-preserving local diffeomorphism
$ F_S:\R\to \R P^1$ and a homomorphism $ f_S:\Z\to \mathrm{PSL}_2(\R)$ that satisfy (\ref{eqqq}).
We can associate the following invariants to $S$: the conjugacy class of
$f_S(1)$ in $\mathrm{PSL}_2(\R)$ and the number
$$\max_{x\in\R,y\in\R P^1}\#( F_S^{-1}(y)\cap [x,x+1]).$$
It can be easily checked that these invariants are sufficient to distinguish any two different structures defined in theorem \ref{circle}.$\diamondsuit$

A {\it projective} structure on the circle is an $(\mathrm{PGL}_2(\R),\R P^1)$-structure. Every M\"obius circle is a projective circle; conversely, one can orient a projective circle so that it will become a M\"obius circle.
\begin{lemma}
\label{projective}
Any two different M\"obius structures on $S^1$ defined in theorem \ref{circle} are not equivalent as projective structures.
\end{lemma}
{\bf Proof of Lemma \ref{projective}.}
We proceed as in the proof of theorem \ref{circle}. Suppose that the projective structure $S$ on $S^1$ is defined by
a couple $( F_S, f_S)$, where $ F_S$ is a local diffeomorphism $\R\to\R P^1$ and $f_S:\Z\mapsto\mathrm{PSL}_2(\R)$ is a homomorphism that satisfy (\ref{eqqq}).
It can be easily checked that the conjugacy class of $ F_S(1)$ in $\mathrm{PGL}_2(\R)$ and the number
$$\max_{x\in\R,y\in\R P^1}\#( F_S^{-1}(y)\cap [x,x+2])$$
depend only on the equivalence class of $S$.
These invariants of projective structures distinguish any two nonequivalent M\"obius structures on $S^1$.$\diamondsuit$

{\bf Remark 2.} Due to Lemma \ref{projective} there exists a bijection $\eu M(S^1)
\leftrightarrow
($projective circles modulo projective diffeomorphisms). Hence, a classification of projective circles can be obtained from theorem \ref{circle}; the classification
given in \cite{kuiper} in not quite correct: if $n>0$, then
$S^1$ provided
with a M\"obius structure of class $P_{n,1}$ and $S^1$ provided
with a M\"obius structure of class $P_{n,-1}$ are not isomorphic as projective circles.

\label{anti}

\begin{figure}\centering
\epsfbox{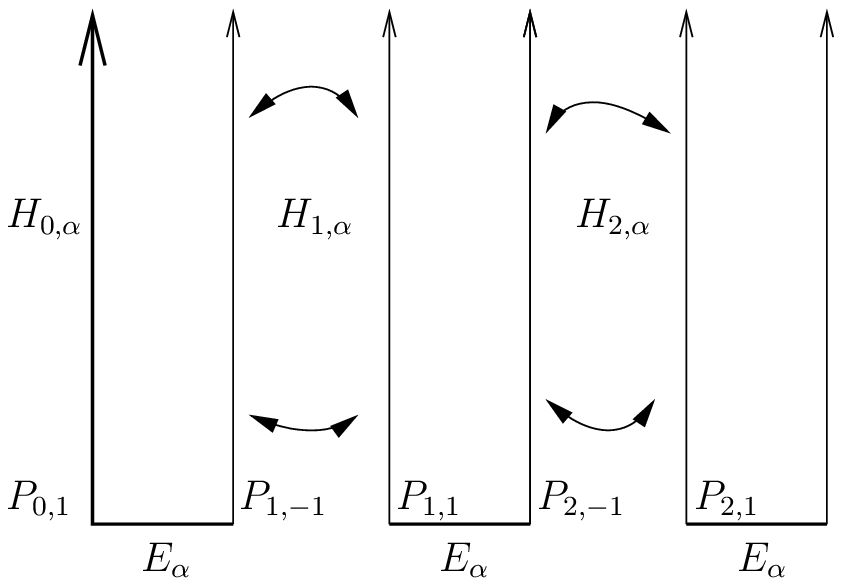}
\caption{$\eu M'(S^1)
$}\label{struct}
\end{figure}

{\bf Remark 3.} The interpretation in terms of developing maps allows us to introduce a topology on $\eu M(S^1)$. However, the resulting topological space
is nasty. The space $\eu M'(S^1)=\eu M(S^1)\setminus\{E_{2\pi k}|\mbox{$k$ is an integer $>0$}\}$
is the non-Hausdorff topological 1-manifold obtained by identifying the thin lines on Figure \ref{struct}; $\eu M'(S^1)$ is open in $\eu M(S^1)$,
and for any $k>0$ we can take the system of sets of the form
$\{E_{2\pi k}\}\cup (\mbox{a neighbourhood of $P_{k,1}$ in $\eu M'(S^1)$})\cup(\mbox{a neighbourhood of $P_{k,-1}$ in $\eu M'(S^1)$})$
as a local neighbourhood basis at $E_{2\pi k}$.


\section{Proof of theorem \ref{nqs}.}\label{sec3}

In order to prove theorem \ref{nqs}, it is enough
to show that any element of $\eu M(S^1)$ is equal to $[S]$, where $S$ is
the M\"obius structure on the boundary of some 2-dimensional M{\"o}bius surface.
Notice that a M\"obius structure on a 2-dimensional manifold
is the same as a $(\mathrm{PSL}_2(\C),\C P^1)$-structure.

\begin{figure}\centering
\epsfbox{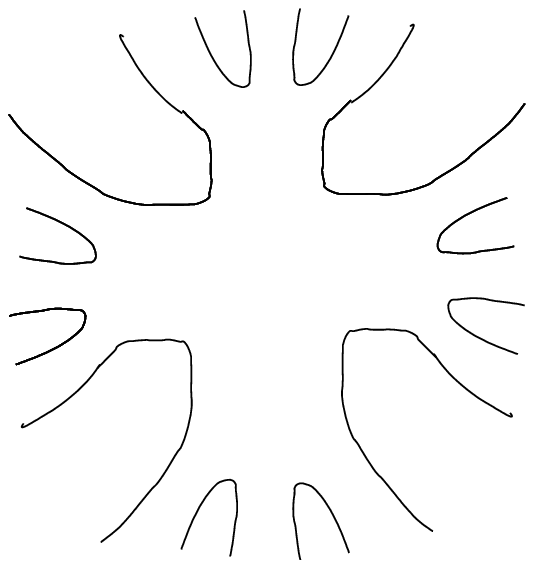}
\caption{}\label{uncov1}
\end{figure}

Let $N$ be the sphere minus three disjoint open disks (``the pants'').
The manifold $\tilde N$ is represented on Figure \ref{uncov1}.

\begin{figure}\centering
\epsfbox{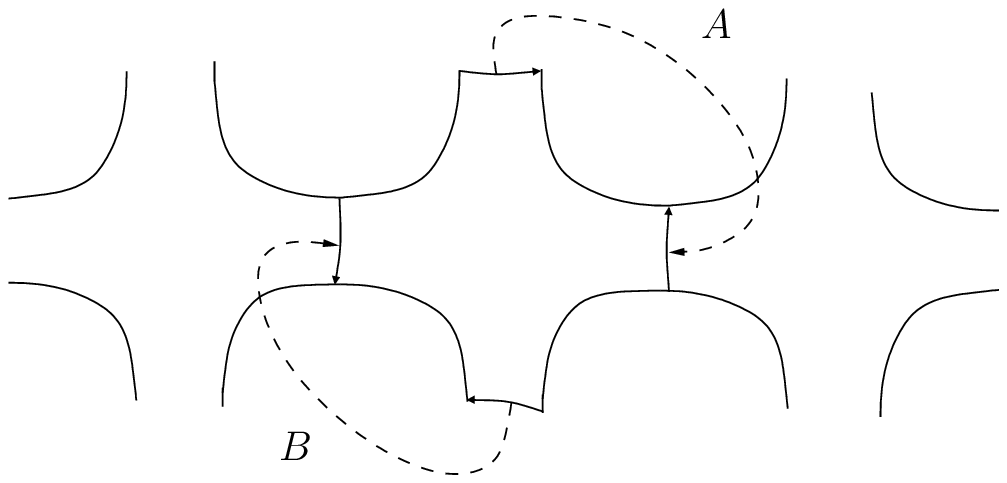}
\caption{}\label{uncov}
\end{figure}

Denote by 
$\mathbf{F}(A,B)$ the free group on the generators $A,B$. 
Suppose that $\mathbf{F}(A,B)$ acts on $\tilde N$ as shown on 
Figure \ref{uncov} (the octagon in the middle
is a fundamental domain; the right and the left octagons are the 
images of the fundamental domain under $A$ and $B$ respectively).

\begin{figure}\centering
\epsfbox{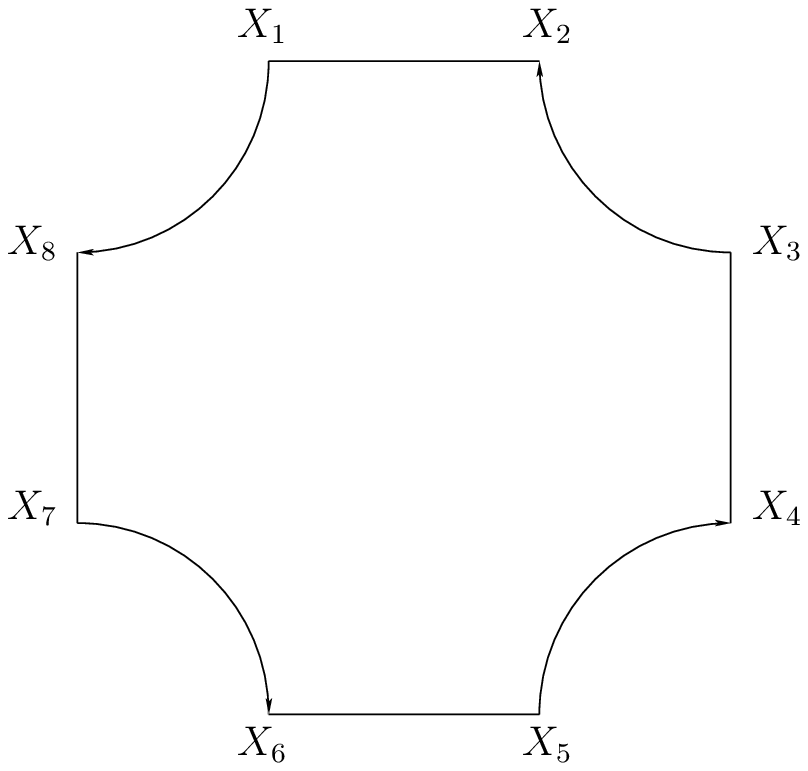}
\caption{}\label{fund}
\end{figure}

We shall equip $N$ (and hence, $\partial N$) with a M\"obius structure by constructing maps
$F:\tilde N\to \C P^1, f:\mathbf{F}(A,B)\to \mathrm{PSL}_2(\C)$ that
satisfy $$F(g\cdot x)=f(g)\cdot x,g\in\mathbf{F}(A,B),x\in\tilde N.$$ This can be done by choosing two
elements $a,b\in \mathrm{PSL}_2(\C)$ and an embedded rectangular octagon 
in $\C$ with vertices $X_1,X_2,\ldots, X_8$ (see Figure \ref{fund}; 
arrows indicate the
orientations of the boundary)
such that the following conditions are satisfied:
\begin{description}\label{cond1}
\item[\rm (C1)] All segments $X_1X_2,\ldots,X_8X_1$ are arcs of circles; 
denote the corresponding circles by $C_{X_1X_2},\ldots,C_{X_8X_1}$.
\item[\rm (C2)] $a(X_1)=X_4, a(X_2)=X_3, b(X_5)=X_8, b(X_6)=X_7$.
\item[\rm (C3)] $a(C_{X_2X_3})=C_{X_2X_3}, a(C_{X_8X_1})=C_{X_4X_5}, 
b(C_{X_6X_7})=C_{X_6X_7}, b(C_{X_4X_5})=C_{X_8X_1}$.
\end{description}
%

\begin{figure}\centering
\epsfbox{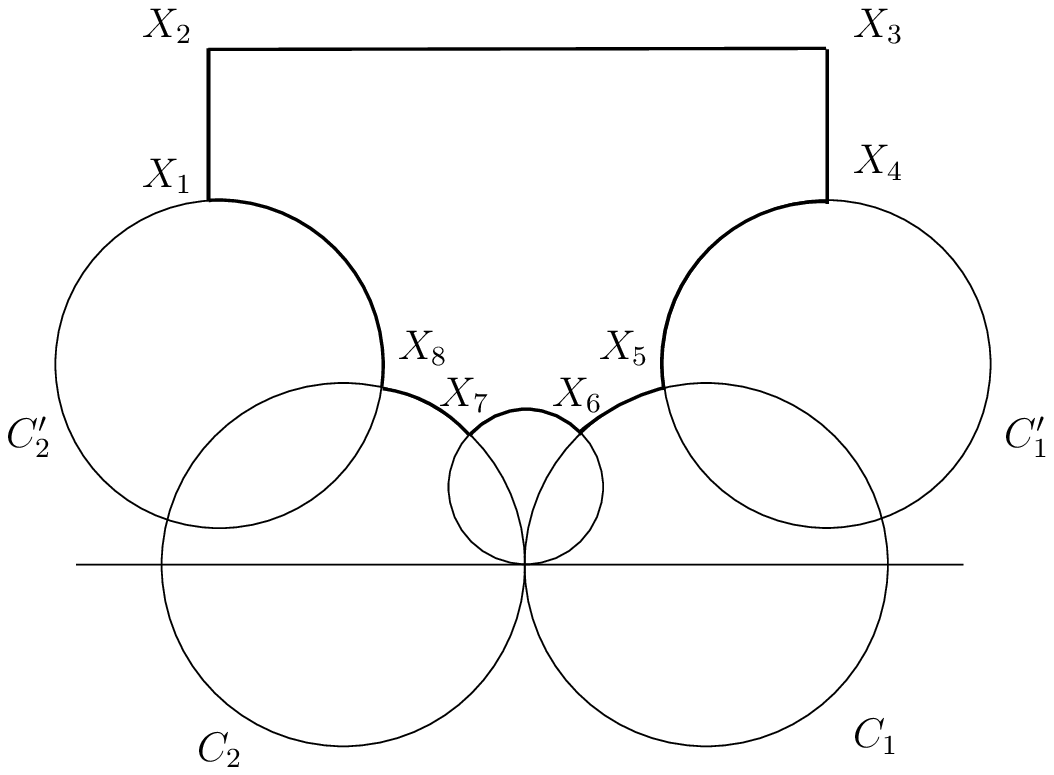}
\caption{}\label{oct}
\end{figure}

Let $x$ be a positive real number. Define the elements 
$a_x, b\in \mathrm{PSL}_2(\C)$ respectively by $z\mapsto z+x, z\mapsto\frac{z}{1-z}$. The transformation $b\circ a_x$ is 
represented by a matrix with trace $2-x$. Choose
the octagon as shown in Figure \ref{oct}. The horizontal line 
at the bottom is the real axis, the circles $C_1$ and $C_2$ on 
Figure \ref{oct} are defined respectively
by the equations $(\mathop{\rm Re}\nolimits z)^2-2\mathop{\rm Re}\nolimits z+(\mathop{\rm Im}\nolimits z)^2=0$ 
and $(\mathop{\rm Re}\nolimits z)^2+2\mathop{\rm Re}\nolimits z+(\mathop{\rm Im}\nolimits z)^2=0$.
The circles $C_1'$ and $C_2'$ are chosen in such a way that 1. $b(C_1')=C_2',a_x(C_2')=C_1'$,
2. $C_1'$ is orthogonal to $C_1$ and depends continuously on $x$, 3.
$C_1'\subset\{z|\mathop{\rm Re}\nolimits z\geq 0\}$, and the imaginary part of
the center of $C'_1$ is $>0$. Note that there are many ways to choose the circles $C_1'$ and $C_2''$ that satisfy these conditions.
The arc $X_6X_7$ is an arc of a sufficiently small circle in the upper half-plane that is tangent to the real axis at $0$.

Note that we have $b(z)=-\bar z$ for $z=1+\mathrm{e}^{it}$, i.e., the 
restriction of $b$ to $C_1$ is the symmetry
with respect to the imaginary axis. This implies easily that 
for any $x>0$, the octagon on Figure \ref{oct} and the 
maps $a=a_x$ and $b$ satisfy the conditions (C1)-(C3).

For any $x>0$ we obtain M{\"o}bius structures on $N$ and all 
components of $\partial N$. The components of $\partial N$ that correspond to $X_2X_3$ and $X_6X_7$
will always have parabolic structures that belong to $P_{0,1}$,
and the equivalence class of the structure on the third component changes as we change $x$. 
Denote this structure by $S_x$; $S_x$ is elliptic for $x<4$, 
parabolic for $x=4$ and hyperbolic for $x>4$.
Denote by $N_1(x)$ the M\"obius surface obtained by gluing together two parabolic components of
$\partial N$ (note that due to  Lemma \ref{projective}, for any M\"obius circle $C$ there exists an anti-M\"obius diffeomorphism $C\to C$).

Let us determine the equivalence class of $S_x$. Let $ F_{S_x}$ be a developing map $\R\to\R P^1$ that corresponds to $S_x$.

Note that $[S_4]=P_{0,1}$, which implies that the image of $ F_{S_4}$ is not the whole $\R P^1$.
Hence, if $D>0$, then we have $ F_{S_x}([0,D])\neq\R P^1$ for any $x$ sufficiently close to 4. This
implies the following lemma:

\begin{lemma}
\label{xx}
For any sufficiently
small $\alpha>0$ there exist $x_1, x_2$ such that the M\"obius structure on $\partial N_1(x_1)$ (respectively, $\partial N_1(x_2)$) belongs to
$H_{0,\alpha}$ (respectively, to $E_{\alpha}$).
\end{lemma}

Let us define the action of $\N$ on $\eu M(S^1)$ as follows: for any $k\in\N, \alpha>0,n\geq 0$ and
$\varepsilon\in\{\pm 1\}$ set
$kH_{n,\alpha}=H_{kn,k\alpha}, kP_{n,\varepsilon}=P_{kn,\varepsilon}, kE_{\alpha}=E_{k\alpha}$.

\begin{lemma}
\label{x}
Let $C'$ and $C''$ be oriented circles. Suppose
that $C''$ is equipped with a M\"obius structure, and denote this structure by $S''$. Let $p:C'\to C''$ be an orientation-preserving $k$-sheeted
covering. Denote by $S'$ the inverse image of $S''$ with respect to $p$. We have $[S']=k[S'']$.
\end{lemma}
$\diamondsuit$
\begin{lemma}
\label{xxxx}
Let $N$ be a compact M\"obius surface with one boundary component and nontrivial fundamental group, and let $k>0$ be an integer.
Denote by $S$ the M\"obius structure on $\partial N$.
There exists
a compact M\"obius surface $N'$ such that $\partial N'$ is a circle, and the M\"obius structure on $\partial N'$ belongs to $k[S]$.
\end{lemma}
{\bf Proof of Lemma \ref{xxxx}.} Let $\mathfrak S_k$ be the symmetric group on $k$ elements. If $k$ is odd, then there exists a homomorphism
$\pi_1(N)\to \mathfrak S_k$ that takes a generator of $\pi_1(\partial N)$ to some cycle of length $k$. Indeed, we can choose
a system $v_1,w_1,\ldots, v_g,w_g$ of free generators of $\pi_1 (N)$ so that $\pi_1 (\partial N)$ will be spanned by
$[v_1,w_1]\cdots [v_g,w_g]$. A required homomorphism can be constructed using the fact that in $\mathfrak{S}_k$ a cycle of length $k$ is a commutator
(e.g., we have
$$(1\ldots k)=\left(1\frac{k+3}{2}\ldots k\right)\left(12\ldots\frac{k+1}{2}\right)=\sigma\left(12\ldots\frac{k+1}{2}\right)^{-1}\sigma^{-1}\left(12\ldots\frac{k+1}{2}\right)$$ for some $\sigma\in\mathfrak{S}_k$).

Hence, for $k$ odd, there exists a surface
$N'$ with one boundary component and a $k$-sheeted covering $N'\to N$.
This, together with Lemmas  \ref{xx} and \ref{x}, allows us to construct for any $\alpha >0$ a M\"obius surface, whose boundary carries 
a M\"obius structure that belongs to $E_{\alpha}$.

Now cut out a small M\"obius disk from $N$. The resulting M\"obius surface has two boundary components; one of these
components has the structure $S$, and the structure on the other one belongs to $E_{2\pi}$.
Note that if $M$ is a compact surface with two boundary components, then a generator
of one of these components can be taken as a free generator of $\pi_1(M)$; hence, for any $k>0$ there exists a $k$-sheeted covering $M'\to M$
that is cyclic over one of the components of $\partial M$. We have already constructed for any integer $l>0$ a M\"obius surface, whose boundary
carries a structure from $E_{2\pi l}$, and the lemma follows.
$\diamondsuit$

\begin{figure}\centering
\epsfbox{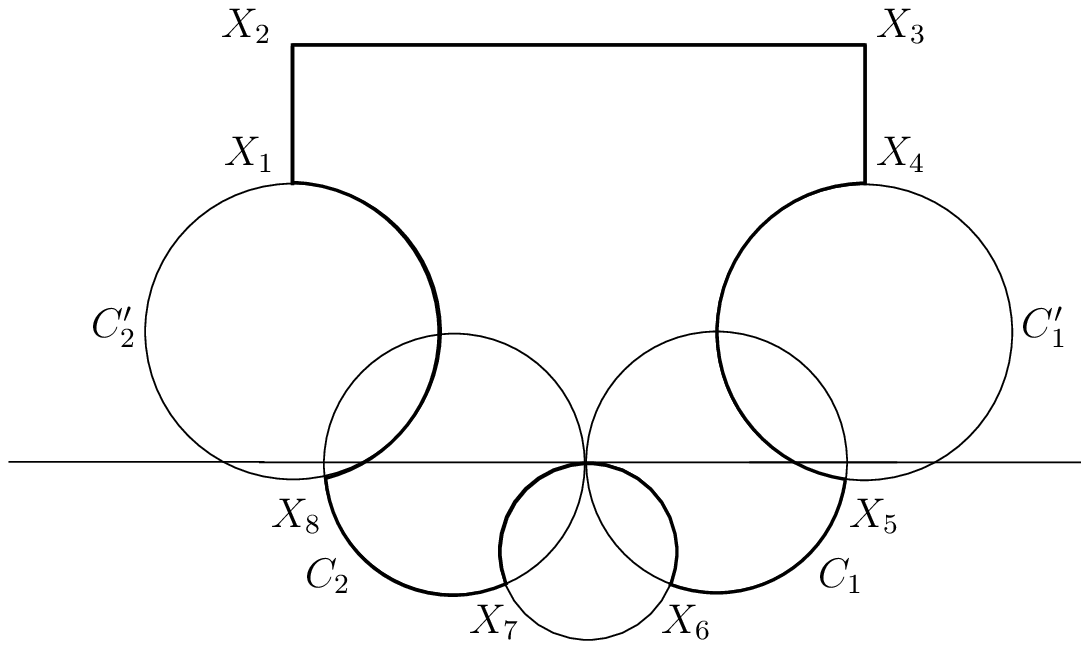}
\caption{}\label{oct2}
\end{figure}

Lemmas \ref{xx}, \ref{x}, \ref{xxxx} imply that any of the classes
$P_{0,1}, H_{0,\alpha}$ or $E_{\alpha}, \alpha>0$ is the equivalence class of the M\"obius structure of the boundary of
some M\"obius surface. In order to complete the proof of theorem \ref{nqs}, we can proceed as follows. Suppose that $x>2$ is a real number and consider
the rectangular octagon represented on Figure \ref{oct2}.

Here the circles $C_1$ and $C_2$ are the same as in the proof of Lemma \ref{xx}. Suppose
that the centers of $C'_1$ and $C'_2$ are the points $\pm \frac x2 +i$, and the radii of these circles are $\frac x2-1$. As above, the octagon
and the maps $a_x$ and $b$ give us for any $x>2$ a M\"obius structure on $N$ (the pants). The induced M\"obius structure
on the component of $\partial N$ that corresponds to
the arc $X_2 X_3$ (respectively, the arc $X_6 X_7$) belongs to the class $P_{0,1}$ (respectively, $P_{1,-1}$). Denote by $S'_x$ the M\"obius structure on the
third component of $\partial N$. The structure $S'_x$ is elliptic for $x<4$; by attaching M\"obius surfaces to the elliptic component of $\partial N$ and to
the $P_{0,1}$-component, we obtain a M\"obius surface, whose boundary has a M\"obius structure of class $P_{1,-1}$, which means that
we can eliminate the $P_{1,-1}$-component as well, i.e., for any $x>2$ there exists a M\"obius surface, whose boundary carries the M\"obius structure
$S'_x$.

Let us determine the equivalence class of $S'_x,x\geq 4$. Note that we can choose a developing map $ F_{S'_x}:\R\to C'_1$ for the structure $S'_x$
so that $ F_{S'_x}$ sends the segment $[0,1/2]$ to the arc $X_5X_4$ and the segment $[1/2, 1]$ to the arc $a_x(X_1 X_8)$. It is easy to check that for $x\geq 4$,
the fixed points of $a_x\circ b$ are exactly the intersection points of $C'_1$ and the real axis.

Hence, for any $x>4$ the set $ F_{S'_x}^{-1}(\mbox{the fixed points of $a_x\circ b$})\cap [0,1)$ consists of two elements, which implies that for any $x>4$ we have
$S'_x\in H_{1,\alpha}$ with $\alpha=2\mathop{\rm arccosh}\nolimits (x/2-1)$.

An analogous argument shows that $S'_4\in P_{1,1}$ (note that $a_4\circ b$ acts on $C'_1$ ``clockwise''). The proof of theorem \ref{nqs} is easily completed
using Lemmas \ref{x} and \ref{xxxx}. $\diamondsuit$
%

\begin{flushright}
{\sc Alexei G. Gorinov\\
IMAPP -- Wiskunde, FNWI\\
Radboud Universiteit Nijmegen\\
The Netherlands}\\
{\tt a.gorinov@math.ru.nl\\
}
\end{flushright}
\end{document}